\documentclass{article}
\usepackage{amsmath,amsfonts,amsthm,amssymb}
\usepackage{graphicx, verbatim}
\newtheorem{theorem}{Theorem}[section]
\newtheorem{corollary}{Corollary}[section]
\newtheorem{lemma}{Lemma}[section]
\title{Uniform Convergence Behavior of the Bernoulli Polynomials}
\author{John Mangual \\
Mathematics Department \\
Princeton University \\
{\tt john.mangual@gmail.com}}
\begin{document}
\maketitle
\begin{abstract} The roots of Bernoulli polynomials, $B_n(z)$, when plotted in the complex plane, accumulate around a peculiar H-shaped curve.  Karl Dilcher proved in 1987 that, on compact subsets of $\mathbb{C}$, the Bernoulli polynomials asymptotically behave like sine or cosine.  Here we establish the asmptotic behavior of $B_n(nz)$, compute the distribution of real roots of Bernoulli polynomials and show that, properly rescaled, the complex roots lie on the curve $e^{- 2\pi  \text{Im}(z)} = 2\pi  e |z|$ or $e^{ 2\pi  \text{Im}(z)}= 2\pi  e |z|$. \end{abstract}

\begin{paragraph}{} This paper came out of the author's Geometry project at Penn State's Mathematics Advanced Study Semesters, Fall 2006.  Adrian Ocneanu showed his geometry class pictures suggesting the roots of Bernoulli polynomials lie on a distinct curve.  See figure 1.  To prove this result, the author learned about the asymptotic behavior of the Bernoulli polynomials.  Observing convergence behavior by plotting the roots of a sequence of polynomials is not new.  Gabor Szeg\"{o} showed in 1928, the roots of $P_n(z) = \sum_{k = 0}^n \frac{x^n}{n!}$ lie on the curve $e^{\text{Re}(z)}=e|z|$.  Jean Dieudonne read Szeg\"{o}'s paper and prove it a different way in 1935.  Also there is a AMS memoir by Karl Dilcher from 1987, \cite{Dil}, where he looks for parabolic zero-free regions in the complex plane.  In this paper we find the exact curve around which the rescaled roots of Bernoulli polynomials cluster.\end{paragraph}
\section{The Bernoulli Polynomials}
\begin{theorem}The roots of $B_n(nz)$ lie on the curve $e^{- 2\pi  \text{Im}(z)} = 2\pi  e |z|$ if $\mathrm{Im}(z) > 0$ and $e^{ 2\pi  \text{Im}(z)}= 2\pi  e |z|$ if $\mathrm{Im}(z) < 0$.
\end{theorem}

\begin{figure}[h]
\centering
\resizebox{!}{50mm}{\includegraphics[height=8.5in]{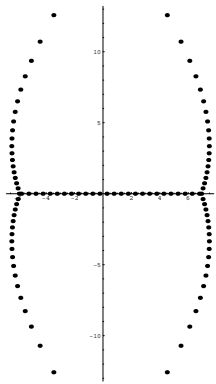}}
\caption{The roots of $B_{100}(z) = 0$ via Mathematica}
\label{moo}
\end{figure} 

\begin{paragraph}{} Let's quickly review some definitions and properties of the Bernoulli numbers and their polynomials.  The Bernoulli numbers are defined by a generating function:
$$ \sum_{k = 0}^\infty B_k \frac{x^k}{k!} = \frac{x}{e^x -1 }$$
There is no simple way to generate the Bernoulli numbers, they are defined as the coefficients of the Taylor expansion of the function above.  The Bernoulli numbers are also coefficients of the Bernoulli polynomials which are also defined by a generating function in two variables:
$$ \sum_{k = 0}^\infty B_k(x) \frac{t^k}{k!} = \frac{te^{xt}}{e^t - 1}$$
One use of the Bernoulli polynomials is to express sums of similar powers:
\begin{equation}(N+1)\sum_{k = 0}^N k^m = B_{m+1}(N+1) - B_{m+1} \end{equation}
In fact Ocneanu's original question as posed in class was to consider the left side of (1), consider it as a polynomial in a complex variable and plot its roots.  This his question was about the sums of like power and not necessarily Bernoulli polynomials.  However, we simplify the problem silghly by not subtracting off terms.\end{paragraph}
\section{History of the Problem}
\begin{paragraph}{}Given a sequence of polynomials, how do we know the roots accumulate on any curve?  Futhermore, how do we compute this curve?  In the October 2005 edition the American Mathemaitcal Monthly I saw this plot:
\end{paragraph}

\begin{figure}[h]
\centering
\resizebox{!}{50mm}{\includegraphics{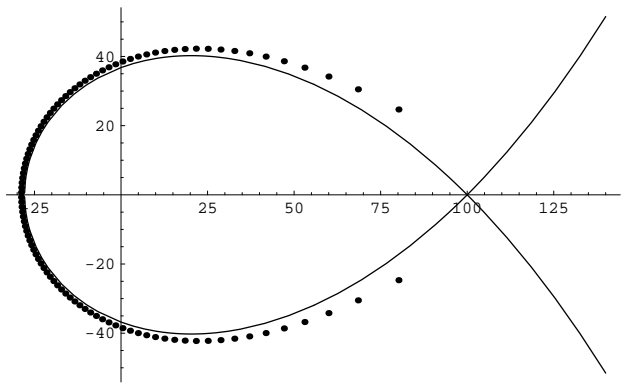}}
\caption{The roots of $\sum_{k = 0}^{100} \frac{(100z)^k}{k!} = 0$ fit the curve similar to $|ze^{z-1}|=1$}
\label{moo}
\end{figure}

\begin{paragraph}{}The paper, \cite{Zem}, was by Stephen Zemyan of Pennsylvania State University, Mont Alto.  He wrote on the roots of the truncated exponential series.  If you cut the taylor series of $e^x$ after $N$ terms you get a polynomial
$$ P_N(x) = \sum_{k = 0}^{N-1} \frac{x^k}{k!}$$
As $N \to \infty$ the roots go to $\infty$ since $e^x$ has no complex roots, but if you rescale by a factor of $N$, the roots of $P_N(Nx)$ approach the curve $|ze^{z-1}| = 1$.  This was proven by Gabor Szeg\"{o} in 1928 appearing in his famous problem book \cite{Sze2} and in a seperate paper \cite{Sze}.  Basically Szeg\"{o}'s result has to do with the gamma function that
$$\sum_{k = 0}^{N-1} \frac{x^k}{k!} = \int_0^x z^N e^z dz$$
However, this work seemed to come out of nowhere and I went in search of more systemattic approach to this problem.  \end{paragraph}
\begin{paragraph}{}Around 1935, Jean Dieudonn\'{e} proved the same result about the rescaled roots of the truncated exponential series. He shows that:
\begin{equation} \frac{n![e^{nz} - P_n(nz)]}{z^n} = \frac{1}{1-z} \left[ 1 + \lambda_n(z)\right] \end{equation} where $\lambda_n(z) \to 0$ uniformly on the $\{z: |z|<1, |z-1|>\epsilon\}$.  Bounding away from $1$ is essential because of the pole in the asymptotics.  In his paper Dieudonn\'{e} says that his solution is a technique which works for general classes of functions.  Having read Dieudonn\'{e}'s proof, the author will try this technique on the Bernoulli polynomials. \end{paragraph}
\begin{paragraph}{} Unlike the exponential problem proved by Szeg\"{o} and Dieudonne, the Bernoulli polynomials are not the truncations of any particular infinite series.  The coeffients seem to ``jump" around and diverge.  There was no obvious way to get around this problem In 1987, Karl Dilcher showed the following formulas using properties of the Riemann Zeta funtion:
\begin{theorem}
As in \cite{Dil2}, let $$T_{2k}(z) = \sum_{j = 0}^k (-1)^j \frac{z^{2j}}{(2j)!}\text{ and }
T_{2k+1}(z) = 
\sum_{j=0}^k (-1)^j\frac{z^{2j+1}}{(2j+1)!}$$ be the truncated sine and cosine series, for $n \in \mathbb{N}$.  Then
 $$\left| (-1)^n \frac{(2\pi)^n}{2n!}B_n\left(z + \frac{1}{2}\right) - T_n(z)\right| < 2^{-n}e^{4\pi |z|}$$
\end{theorem}

\begin{corollary} On compact subsets of $\mathbb{C}$ we get uniform convergence:
$$ \lim_{n \to \infty} \frac{(2\pi)^n}{2n!}B_n(z) \to   \left\{ \begin{array}{cc}
 \cos\left(2\pi z\right) & \text{if }n\text{ is even}\\ \\
 \sin\left(2\pi z\right) & \text{if }n\text{ is odd}\end{array}\right.$$
\end{corollary}
This showed the Bernoulli polynomials were intimately related to the basic trigonometric functions and gave a quantative error between the limits of the respective truncations.  However, if we substitute $nz$ as an argument instead of $z$, we get uniform convergence on a small disc where $e^{4\pi|z|}/2 < 1$ or $|z| < \ln 2 /4\pi = 0.055 < 0.159 = 1/2 \pi$.  This means we cannot show these functions are unformly close on the full disc $\frac{1}{2\pi}\mathbb{D}$ on which the roots lie.  However, we do get a compelling picture
\begin{figure}[h]
\centering
\resizebox{!}{50mm}{\includegraphics{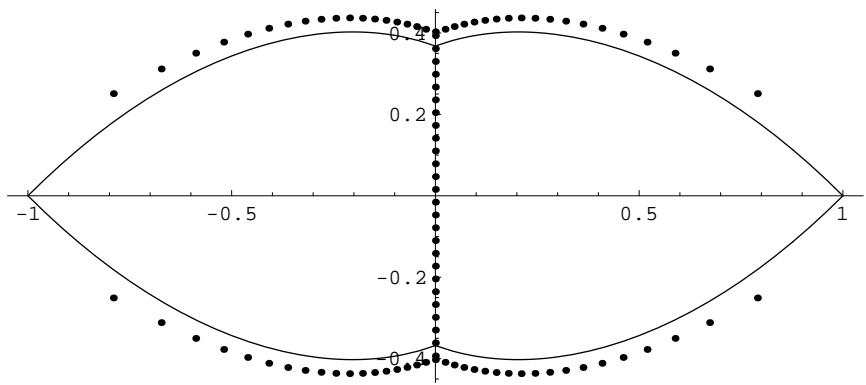}}
\caption{The roots of the Bernoulli polynomials compared to the locus of roots of truncated cosine.}
\label{moo}
\end{figure}
\end{paragraph}
\section{The Real Case}
\begin{paragraph}{}Now we will show the real roots of $B_n(nz)$ are uniformly spaced on an interval $[-1/2\pi n e, 1/2\pi n e]$.  Since Karl Dilcher related Bernoulli polynomials to trigonometric series and exponential functions, Fourier series may be the appropriate tool.  Let's recall a result of Hurwitz:\end{paragraph}
\begin{lemma} The Fourier series of the Bernoulli polynomials are given by:
 \begin{equation} B_n(\{x\}) = - \frac{n!}{(2\pi i)^n} \sum_{|k|>0} \frac{e^{2\pi i  k x}}{k^n}\end{equation} 
Notice this Fourier series has period 1 so we use the fractional part sign $\{x\}$. \end{lemma}
 \begin{proof}
The coefficients can be found all at once using the generating function:
$$\int_0^1 \frac{te^{xt}}{e^t - 1} e^{-2\pi i n x} dx = \frac{t}{t - 2\pi i n } = -\sum_{k = 1}^\infty \left(\frac{t}{2\pi i n}\right)^k$$ For each individual $k$ we get a Fourier series like the one above.
\end{proof}

\begin{paragraph}{} What are the asymptotics of this periodized version of the Bernoulli polynomials? \end{paragraph}
\begin{lemma} The Bernoulli polynomials $B_n(\{x\})$ converge uniformly to $$-\frac{(2\pi i )^n}{2n!}\times \frac{e^{2\pi i x} + e^{-2\pi i x}}{2}$$
in the unit disc $\mathbb{D}$. \end{lemma}

\begin{proof} Consider the Fourier series from the previous lemma:
$$ B_n(\{x\}) = - \frac{n!}{(2\pi i)^n} \sum_{|k|>0} \frac{e^{2\pi i  k x}}{k^n}$$
Every term except $k = \pm 1$ goes to zero as $n$ goes to infinity.  This gives us pointwise convergence.  Is this convergence uniform for  $x \in \mathbb{D}$? Yes
$$ \left| \sum_{|k|\geq 2} \frac{e^{2\pi i  k x}}{k^n}\right| < 2 \sum_{k\geq 2} \frac{1}{k^n} < C \int_2^\infty \frac{dx}{x^n} = \frac{C'}{2^n}$$
\end{proof}
\begin{paragraph}{} Now we have estimated the periodized version of the Bernoulli polynomials $B_n(\{z\})$ but how do we estimate $B_n(mx) - B_n(\{mx\})$? \end{paragraph}
\begin{lemma} As $n$ approaches infinity $B_n(mx) - B_n(\{mx\})$ approaches
$$  nm^{n-1} \int_0^x mt^{n-1} dt $$
This convergence is uniform in the unit interval $[0, 1-\epsilon]$ for any $\epsilon > 0$. \end{lemma} 
\begin{proof}Let's do one at a time: $B_n(x+1) - B_n(x) =  nx^{n-1}$.
Then 
\begin{eqnarray*} B_n(mx) - B_n(\{mx\}) &=& \sum_{k = 0}^{\lfloor nx\rfloor - 1} B_n(\{ mx \} + k + 1) - B_n( \{mx\} + k\} \\
&=&   \sum_{k = 0}^{\lfloor nx\rfloor - 1} n\left( \{mx\} + k  \right)^{n-1} \\
&=&  nm^{n-1} \times \frac{1}{n} \sum_{k = 0}^{\lfloor nx\rfloor - 1} m\left( \frac{\{mx\} + k}{m}  \right)^{n-1} \\
&\approx & nm^{n-1} \times \int_0^x mt^{n-1} dt\end{eqnarray*}
For $x \in [0,1]$, does the Riemann sum above convergence uniformly in $\mathbb{D}$ to the following integral?
$$ \int_0^x mt^{n-1} dt = x^n$$
Since $x^m$ is monotone in $x$ on $[0,1]$, the difference between the upper and lower Riemann sums is just the difference between the first term of the upper sum and the last term of the lower sum.  The Riemann sum above is the lower sum.  The upper sum is defined between $k=1$ and $k = \lfloor n x \rfloor $.  Their difference is just $x^{n-1} -  [\{mx\}/m]^{n-1} < x^{n-1}$.  This means convergence will be uniform for $x \in [0, 1-\epsilon]$ for $\epsilon > 0$.  
\end{proof}
\begin{paragraph}{} How let's take another look at (3):
$$B_n(\{x\}) = - \frac{n!}{(2\pi i )^n}\sum_{k = -\infty}^\infty \frac{e^{2\pi i n kx}}{k^n}$$
We need an estimate of Stirling's formula with errors.  By considering the integral
$$ f(t) = \int_0^\infty e^{t(\ln x - x)} dx $$
it is possible using Laplace's method to get Stirling's esimate:
$$ n! = (n/e)^n \sqrt{2\pi n } \left[ 1 + \frac{1}{12n} + \frac{\mathcal{O}(1)}{n^2}\right]$$
Thus
\begin{eqnarray*} B_n(\{x\}) &=& - \frac{(n/e)^n \sqrt{2\pi n } \left[ 1 + \frac{1}{12n} + \frac{\mathcal{O}(1)}{n^2}\right]}{(2\pi i )^n}\sum_{k = -\infty}^\infty \frac{e^{2\pi i n kx}}{k^n} \\ & \approx & - \frac{[n\hspace{0.1em}\mathcal{O}(n^{1/2n})]^n}{(2\pi i e)^n}(e^{2\pi i nx} - e^{-2\pi i n x} + \mathcal{O}(2^{-n})) \end{eqnarray*}
Asymptotically the roots of $B_n(nx) = B_n(\{nx\}) + (B_n(nx) - B_n(\{nx\}))$ are the roots of $e^{2\pi i n x} - e^{-2\pi i n x} = 0$.  This is basically the stategy we will use in passing to the complex roots.
\end{paragraph}
\section{The Complex Case}
\begin{paragraph}{} The previous section's result can be found in \cite{Ves} from 1999.  In this section, while the proof is new, the result was anticipated by Karl Dilcher.  Let's start writing the simple generating function relation
$$ \frac{te^{xt}}{e^t - 1} \times e^{yit} = \frac{te^{(x + yi)t}}{e^t - 1}$$
and rewriting it in terms of Bernoulli polynomials:
$$ B_n(x + yi) = \sum_{k = 0}^n \binom{n}{k}B_k(x)(yi)^{n-k}$$
Let's define a functional:
$$ \sum_{k = 0}^\infty F_n(x)\frac{t^n}{n!} \to e^{yit} \sum_{k = 0}^\infty F_n(x)\frac{t^n}{n!} \hspace{0.5cm}\text{ or }\hspace{0.5cm}F_n(x) \to \sum_{k = 0}^n \binom{n}{k}F_k(x)(yi)^{n-k}$$
This functional is linear in generating functions $F(x,t)$ or equivalently, sequences $F_n$.  In particular, we evaluate this functional for $F_n(x) = B_n(\{nx\})$ and $F_n(x) = B_n(nx) - B_n(\{nx\})$.\end{paragraph}
\begin{lemma}As $n$ approaches infinity,  $B_n(nx+ niy) - B_n(\{nx\}+iny)$ approaches
$$ n^n \int_0^1 n(x + yi)^{n-1}dx $$ \end{lemma}
\begin{proof}
 \begin{eqnarray*} \sum_{k = 0}^n \binom{n}{k} (B_k(nx) - B_k(\{ nx\} ))(yi)^{n-k} &=& n^n \sum_{k = 0}^n \binom{n}{k} 
 \times \frac{1}{n} \sum_{k = 0}^{\lfloor nx\rfloor - 1} n\left( \frac{\{ nx\} + k}{n}  \right)^{n-1}
 (yi)^{n-k}\\
 &=& n^n \times  \frac{1}{n} \sum_{k = 0}^{\lfloor nx\rfloor - 1} n\left( \frac{\{ nx\} + k}{n}  + yi \right)^{n-1}
 \end{eqnarray*}
 This cooresponds to a Riemann integral:
 $$\int_0^1 n(x + yi)^{n-1}dx$$
 and since this function is monotone in $x$ we get the same error bound as before:
 $$ n(x + yi)^{n-1} - n\left(\frac{\{nx\}}{n} + yi\right)^{n-1}$$
 As long as $|x + yi|< 1- \epsilon$ we get uniform convergence for $\epsilon > 0$.
 \end{proof}
 Now let's give our final major estimates of the paper:
\begin{lemma}As $n$ approaches infinity $B_n(\{nx\} + iny)$ approaches
$$- \frac{[n\hspace{0.1em}\mathcal{O}(n^{1/2n}(1 + 1/n))]^n}{(2\pi i e)^n} \left[ P_n(2 \pi ny) + P_n(-2\pi n y)  \right] $$\end{lemma}
\begin{proof} Now we transform $B_n(\{nx\})$.  We already know its Fourier series expansion for real arguments.  Let's use the convolution formula above to get add an imaginary value $iny$ :
\begin{eqnarray*}B_n(\{nx\} + iny ) &=& \sum_{m=0}^n \binom{n}{m}B_m(\{nx\}) (iny)^{n-m}\\ &=& -
n! \sum_{m=0}^n  \frac{1}{(2\pi i)^m}  \sum_{k = -\infty}^\infty  k^{-m}e^{2\pi i k nx}\frac{ y^{n-m}}{(n-m)!}\\
&=& -
\frac{n!}{(2\pi i )^n} \sum_{m=0}^n  \frac{ (2\pi i y)^{n-m}}{(n-m)!}  \sum_{k = -\infty}^\infty  k^{-m}e^{2\pi i k n x}\end{eqnarray*}
Therefore we get the following estimate for $B_n({nx} + iny)$.  We pull out the $|k|=1$ terms in the Fourier series, which do not decay with $m$.  In the third line, we use our uniform bound for the $|k|>2$ tail in Lemma 3.2 to give line four.  
 \begin{eqnarray*}B_n(\{nx\} + iny) &=& -
\frac{n!}{(2\pi i )^n} \sum_{m=0}^n  \frac{ (-2\pi n y)^{n-m}}{(n-m)!}  \sum_{k = -\infty}^\infty  k^{-m}e^{2\pi i k n x}\\ &=& -
\frac{n!}{(2\pi i )^n} \sum_{m=0}^n  \frac{ (-2\pi n y)^{n-m}}{(n-m)!} \left[ e^{2\pi i n x} + (-1)^ne^{-2\pi i n x} +  \sum_{|k|\geq 2}  k^{-m}e^{2\pi i k n x} \right] \\
&=& -
\frac{n!}{(2\pi i )^n} \sum_{m=0}^n  \frac{ (-2\pi n y)^{n-m}}{(n-m)!} \left[ e^{2\pi i n x} + (-1)^ne^{-2\pi i n x} +  \mathcal{O}(2^{-m}) \right] \\
&=& -
\frac{n!}{(2\pi i )^n} \left[ P_n(2 \pi ny) + P_n(-2\pi n y) + \frac{1}{2^n}\mathcal{O}(P_n(\pm 4\pi ny)) \right] \end{eqnarray*}
Having summed all the errors we get a worrisome error which might grow as $\frac{1}{2^n}\mathcal{O}(P_n(\pm 4\pi ny))$ which might grow faster than the dominant terms!   Indeed our proof would collapse if we estimated this as $\frac{1}{2^n}\mathcal{O}(e^{\pm 4\pi ny})$, but it so happens that $\frac{1}{2^n}P_n(4\pi n |y|) < P_n(2\pi n |y|)$.  Let's double-check this:
$$\frac{1}{2^n}P_n(4\pi n |y|) = \frac{1}{2^n}\sum_{k = 0}^{n-1} \frac{(4\pi n |y|)^k}{k!}
< \sum_{k = 0}^{n-1} \frac{(2\pi n |y|)^k}{k!} = P_n(2\pi n |y|)$$
\end{proof}
Now we use Dieudonne's estimate to say that within the unit circle our Taylor series converges to $P_n(2\pi i nz) \to e^{2\pi i nz}$ uniformly inside the disc $\mathbb{D}/2\pi$.
$$ n![e^{2\pi i nz} - P_n(2\pi inz)] = \frac{(2\pi i z)^n}{1-2\pi i z} \left[ 1 + \lambda_n(z)\right] $$
This is to say the $P_n(2\pi i n z)$ converges to the exponential uniformly on the disc $\frac{1}{2\pi}\mathbb{D}$.  Also note that $[x^n(1 + \mathcal{O}(1)]^{1/n} = x[1 + \mathcal{O}(1/n)]$.
\begin{equation}B_n(\{nx\} + iny) =  -
\frac{n!}{(2\pi i )^n} \left[ e^{2 \pi ny} + e^{-2\pi n y}\right] [1 + \mathcal{O}(1) ]
=  - \frac{[n\hspace{0.1em}\mathcal{O}(n^{1/2n}(1 + 1/n))]^n}{(2\pi i e)^n} \left[ e^{2 \pi ny} + e^{-2\pi n y}  \right]\end{equation}
and we use Stirling's formula to get the above estimate.

\begin{paragraph}{} Earlier we showed:
$$B_n(nx+ niy) - B_n(\{nx\}+iny) =  [n(x + iy)]^n\left[1 + \left(\frac{x}{x + iy}\right)^n\right] $$
This tells us this equation we are looking for is:
$$2\pi ie z = \left\{ \begin{array}{cc} e^{2\pi n y} & \text{if }y > 0 \\
e^{-2\pi n y} & \text{if }y < 0 \end{array} \right. $$
 \end{paragraph}

\section{Acknowledgement}
\begin{paragraph}{} Many thanks to Penn State especially Adrian Ocneanu, Sergei Tabachnikov and Alberto Bressan.  Also thanks to Karl Dilcher himself who wrote an entire memoir, \cite{Dil2} on the roots of Bernoulli polynomials in 1987.  He was responsive to my e-mails.  Also thanks to Stepen Zemyan who showed me how to plot roots of polynomials.
\end{paragraph}

\end{document}